\documentclass[11pt]{amsart}
\usepackage{amssymb}
\newtheorem{theorem}{Theorem}[section]
\newtheorem{proposition}[theorem]{Proposition}
\newtheorem{lemma}[theorem]{Lemma}
\theoremstyle{definition}
\newtheorem{definition}[theorem]{Definition}
\newtheorem{notation}[theorem]{Notation}
\newcommand{\PP}{\mathbb P}
\newcommand{\Aff}{\mathbb A}
\newcommand{\CC}{\mathbb C}
\newcommand{\OO}{\mathcal O}
\newcommand{\Dom}{\operatorname{Dom}}
\newtheorem{mainresult}{Theorem}

\newcommand{\HypersurfaceStatement}{There is a dense Zariski-open subset of
$|\OO_{\PP^2\times\PP^2}(2,2)|$ whose members are smooth algebraically
elliptic Fano threefolds. A very general member is also stably irrational.}
\newcommand{\CoverStatement}{Assume that $S_B$ is smooth and that
$\Delta_f$ and $\Delta_g$ are smooth plane sextics. Then $X_B$ is
algebraically elliptic. These hypotheses hold for a general section $B$.
For a very general section $B$, the threefold $X_B$ is also stably irrational.}
\title[Ellipticity and two conic bundles]{Algebraic ellipticity of threefolds with two conic bundle structures}
\author{Mikhail Zaidenberg}
\address{Univ. Grenoble Alpes, CNRS, IF,
38000 Grenoble, France}
\email{mikhail.zaidenberg@gmail.com}
\date{}
\expandafter\providecommand\csname subjclassname@2020\endcsname{2020 Mathematics Subject Classification}
\begin{document}
\begin{abstract}
We prove algebraic ellipticity for a general hypersurface of bidegree $(2,2)$
in $\PP^2\times\PP^2$ and for a double cover of the incidence threefold
branched along a smooth divisor of bidegree $(2,2)$ when the two conic
bundle discriminants are smooth sextics. These are the two cases of
Mori--Mukai family~$2.6$. A very general member of each of the two families
is stably irrational. We also prove rationality and algebraic ellipticity for smooth hypersurfaces with two conic bundle projections in products of smooth rational surfaces whenever at least one surface is different from $\PP^2$.
\end{abstract}
\maketitle

\section*{Introduction}
We work over $\CC$. A smooth variety $X$ is said to be algebraically elliptic 
in the sense of Gromov if there are an algebraic vector bundle $E\to X$ and a morphism $s:E\to X$ 
such that $s(0_x)=x$ and the differential of $s|_{E_x}$ at $0_x$ surjects onto $T_xX$ for every $x$, 
see \cite{Gro}. Algebraic ellipticity implies many additional properties, see e.g. \cite{For} and \cite{Zai}. 
As Kaliman observed \cite{KalBan}, Banecki's proof for smooth projective rational varieties 
\cite[Corollary~1.9]{Ban} gives algebraic ellipticity for every smooth complete rational variety. 
For recent examples of irrational and stably irrational algebraically elliptic threefolds 
and fourfolds see \cite{DKZ26a},  \cite{DKZ26b}, \cite{DKZ26c}, and  \cite{KZcubic}.

We consider the Fano threefolds of Mori--Mukai family~$2.6$. They have
Picard number two, index one, anticanonical degree twelve, and genus seven.
In case~$(a)$ the two conic bundle projections together embed the threefold
as a hypersurface of bidegree $(2,2)$ in a product of two planes
\cite[Section~4.1]{Sharma}. B\"ohning and Graf von Bothmer proved that a
very general such hypersurface is stably irrational \cite{BB}.

Case~$(b)$ is a specialization of case~$(a)$ \cite[\S7.2]{HTF}.
Let $P=\PP^2$, let $P^\vee$ be its dual plane, and let
$F=\{(p,L)\in P\times P^\vee\mid p\in L\}$ be the incidence threefold.
For $B\in H^0(F,\OO_F(2,2))$ let $S_B$ be its zero divisor and
$\pi:X_B\to F$ the double cover defined by $w^2=B$ in the total space of
$\OO_F(1,1)$. Composing $\pi$ with the two projections gives morphisms
$f:X_B\to P$ and $g:X_B\to P^\vee$ with plane conics as fibers.
Denote their discriminants by $\Delta_f$ and $\Delta_g$.
Hassett and Tschinkel proved that a very general such double cover
is stably irrational \cite[\S7.2, following Lemma~7]{HTF}.

By a conic bundle we mean a flat projective morphism whose general fiber is a smooth conic and whose geometric fibers are all isomorphic, as schemes, to plane conics. In particular, a singular geometric fiber is either a pair of distinct lines or a double line.

\begin{mainresult}\label{main:unified}
Let $S_1$ and $S_2$ be smooth rational projective surfaces. In every complete linear system on $S_1\times S_2$ containing a smooth hypersurface whose projections are conic bundles, a general such hypersurface is algebraically elliptic. More precisely, the following statements hold for such hypersurfaces and for double covers of the incidence threefold.
\begin{enumerate}
\item If at least one of the two surfaces is different from $\PP^2$, then every such hypersurface $X$ is rational and algebraically elliptic. The same conclusion holds when $S_1=S_2=\PP^2$ and the bidegree of $X$ is either $(1,n)$ or $(n,1)$.
\item \HypersurfaceStatement
\item \CoverStatement
\end{enumerate}
\end{mainresult}

We use the local criterion for algebraic ellipticity \cite{KZ} in the
form of \cite[Lemma~2.5]{KZcubic}. A local rank-one spray at $x\in X$
is given, after trivializing its line bundle, by a morphism
$s:U\times\Aff^1\to X$ on a Zariski neighborhood $U$ of $x$, with
$s(u,0)=u$. Its velocity at $x$ is the tangent vector to the curve
$t\mapsto s(x,t)$ at $t=0$. By \cite[Lemma~2.5]{KZcubic}, it suffices to construct at each $x$ a
family of such sprays whose velocities span $T_xX$.

In the two Fano cases we write $C_b=f^{-1}(b)$ and $D_p=g^{-1}(p)$.
Suppose first that $x$ is smooth on $C_{f(x)}$. A spray along this
conic supplies a nonzero velocity in $\ker df_x$. We use birational
reflections preserving $f$ to transport two further sprays to $x$
from suitable points of the same conic. We choose these points so
that the two transported velocities project to distinct tangent
directions in the base of $f$. The resulting three velocities span
$T_xX$. Our use of reflections follows
\cite{KZcubic,DKZ26a,DKZ26b,DKZ26c}.

For the hypersurface in Theorem~\ref{main:unified}(2), smoothness of $X$
ensures that at each point at least one of the two conics is smooth there.
Thus the preceding construction, with the projections interchanged when
necessary, covers every point. For the incidence double cover in
Theorem~\ref{main:unified}(3), a point can be singular on both conics.
At such a point $x$, a pencil of lines through $f(x)$ supplies a general
line $\ell$. We choose a second point $d\in\ell\cap\Delta_f$. The pencil of lines through $d$ gives a family of surfaces $f^{-1}(m)$, with $f^{-1}(\ell)$ as its distinguished member. Two successive contractions give a family of degree-four
del Pezzo surfaces. On this family we construct sprays whose images
avoid the section over which the contraction fails to be an isomorphism.
They lift to give two velocities spanning $\ker df_x$. Interchanging
the projections gives two more spanning $\ker dg_x$, and these two
planes together span $T_xX$.

In Section~\ref{sec:rational-surfaces} we prove the rationality assertion of Theorem~\ref{main:unified}(1) by showing that one of the projections admits a rational section. Algebraic ellipticity then follows from \cite[Corollary~1.9]{Ban}.

\smallskip

\noindent\textbf{AI disclosure.}
This work grew out of discussions between the author and ChatGPT Astra~6.2.
The author contributed ideas for the proof, and Astra~6.2 developed rigorous
arguments from those ideas. The author and Astra~6.2 subsequently revised
the manuscript together. The author takes responsibility for its mathematical content.

\section{Hypersurfaces of bidegree $(2,2)$ in a product of planes}

\begin{theorem}[Theorem~\ref{main:unified}(2)]\label{thm:biquadratic}
\HypersurfaceStatement
\end{theorem}

\subsection{The two conic bundles} \label{sec:2}

Write $F(a,z)=0$ for the equation of $X$ in the product of projective planes $P_a\times P_z$. 
Thus $F$ is a bihomogeneous biquadratic form in $3+3$ variables.
The projections $f:X\to P_a$ and $g:X\to P_z$ define two conic bundles on $X$.  
For $b\in P_a$ and $q\in P_z$, set
\begin{equation}\label{eq:conics}
 Q_b(z)=F(b,z),\quad C_b=\{Q_b=0\}\subset P_z,\quad
 D_q=\{F(a,q)=0\}\subset P_a.
\end{equation}
We identify $C_b$ with $\{b\}\times C_b$ in $X$, and $D_q$ with $D_q\times\{q\}$. The discriminant $\Delta_f\subset P_a$ is the locus of $b\in P_a$ for which $C_b$ is singular. We define $\Delta_g\subset P_z$ in the same way, and we choose affine coordinates $(u,v)$ in $P_a$ centered at $b$ and a local trivialization of $\OO_{P_a}(2)$. Viewing $F(a,z)$ as a quadratic form in $z$, we differentiate its coefficients with respect to $u$ and $v$ and evaluate at $a=b$. The resulting quadratic forms in $z$ are denoted by $Q_u$ and $Q_v$, respectively. In the chosen trivialization, the expansion to first order at $b$ is
\begin{equation}\label{eq:expansion}
F(a,z)\equiv Q_b(z)+u\,Q_u(z)+v\,Q_v(z)
\pmod{(u,v)^2}.
\end{equation}
Changing the trivialization multiplies their restrictions to $C_b$ by the same nonzero scalar, since the extra terms are multiples of $Q_b$. Changing the base coordinates makes an invertible linear change in the pair.

\begin{lemma}\label{singular-fibers}
For a general form $F$ of bidegree $(2,2)$, the hypersurface $X$ is smooth, and every singular fiber of either projection is the union of two distinct lines.
\end{lemma}
\begin{proof}
Smoothness follows from Bertini's theorem. Set 
\[V=H^0(P_a\times P_z,\OO(2,2)),\quad\text{where}\quad\OO(2,2)=\OO_{P_a\times P_z}(2,2),\] and let $A_F(b)$ be the symmetric matrix of $Q_b$. At a fixed $b$, the three quadrics $Q_b,Q_u,Q_v$ can be prescribed independently, since they are given by the coefficients of $1,u,v$ in the first-order expansion \eqref{eq:expansion} of $F$. Indeed, in homogeneous coordinates $[a_0:a_1:a_2]$ with $b=[1:0:0]$, the coefficients of $a_0^2,a_0a_1,a_0a_2$ in $F$ can be chosen as arbitrary quadratic forms in $z$. In particular, $A_F(b)$ ranges over all symmetric $3\times3$ matrices. Those of rank at most one have codimension three. The closed incidence variety of pairs $(F,b)$ with $\operatorname{rank}A_F(b)\leq1$ therefore has dimension $\dim V-1$, since $b$ varies in a plane. Its image in $V$ is closed by projectivity of $P_a$. Outside this image every singular fiber has rank two. Interchanging the factors gives the same conclusion for the second conic bundle on a common dense open subset of $|\OO(2,2)|$.
\end{proof}

\begin{lemma}\label{discriminants}
For a general $F$, the discriminants of both projections are smooth irreducible plane sextics.
\end{lemma}
\begin{proof}
We retain the notation of Lemma~\ref{singular-fibers}. The equation of $\Delta_f$ is $\det A_F(b)=0$. Its degree is six, since the matrix entries are quadratic forms in $b$, and it is nonzero for general $F$, since we may prescribe a smooth conic over one fixed point. At any $b\in\Delta_f$, Lemma~\ref{singular-fibers} allows us to choose coordinates with $A_F(b)=\operatorname{diag}(1,1,0)$. The node of the fiber is then represented by $z_0=(0,0,1)$. In the affine chart $z=(s,t,1)$, the equation of $C_b$ is $s^2+t^2=0$, so both partial derivatives with respect to $s,t$ vanish at the node. Since $X$ is smooth there, its local equation has a nonzero partial derivative with respect to $u$ or $v$. By~\eqref{eq:expansion}, these two partial derivatives are $Q_u(z_0)$ and $Q_v(z_0)$. They are also the corresponding partial derivatives of $\det A_F$, since the only nonzero cofactor of $A_F(b)=\operatorname{diag}(1,1,0)$ is its $(3,3)$-cofactor, which equals one. Thus $\Delta_f$ is smooth. It is irreducible because any two plane curves of positive degree meet. The same proof applies to $\Delta_g$.
\end{proof}

\begin{lemma}\label{derivative-quadrics}
For a general $F$, the restrictions of $Q_u,Q_v$ to each smooth fiber $C_b$, 
and to every component of each reducible fiber $C_b$, are linearly independent. 
The same holds for the second conic bundle.
\end{lemma}
\begin{proof}
Recall that $Q_b(z)=F(b,z)$. Written in the six monomials $z_i z_j$, $0\leq i\leq j\leq 2$, the coefficients of $Q_b$ are quadratic forms in the homogeneous coordinates $b_0,b_1,b_2$ of $b$. Each is therefore a linear combination of the six monomials $b_k b_\ell$, $0\leq k\leq\ell\leq 2$. For general $F$, the resulting $6\times6$ coefficient matrix is invertible, so the map $P_a\to\PP^5$, $b\mapsto[Q_b]$, is the quadratic Veronese embedding followed by a projective automorphism of $\PP^5$. Its differential has rank two at $b$, and by~\eqref{eq:expansion} its values on the two coordinate tangent vectors are represented by $Q_u$ and $Q_v$ modulo the span of $Q_b$. Thus $Q_u,Q_v$ are independent modulo $Q_b$, which proves the assertion on a smooth $C_b$, since every quadric vanishing on $C_b$ is a multiple of $Q_b$.

At fixed $b$, the condition that $Q_b$ have rank two has codimension one, and choosing one of its two component lines, say $M$, adds no dimension. By the independent prescription of quadratic forms as in the proof of Lemma~\ref{singular-fibers}, $Q_u|_M,Q_v|_M$ can be arbitrary binary quadrics. Dependent pairs have codimension two, since their coefficient matrices are the $3\times2$ matrices of rank at most one. Dependence on some component therefore imposes three conditions at fixed $b$. Allowing $b$ to vary adds only two dimensions, so a general $F$ avoids this condition for every $b$. Interchanging the factors proves the assertion for the second conic bundle.
\end{proof}

\subsection{Birational reflections of a conic bundle}\label{sec:3}

We keep the notation of Subsection~\ref{sec:2}, and we fix a biquadratic form $F$ 
satisfying the conclusions of Lemmas~\ref{singular-fibers}--\ref{derivative-quadrics}. 
Let $B_b$ be the symmetric bilinear form with $B_b(z,z)=Q_b(z)$.

\begin{notation}\label{reflection}
Fix $c\in P_z$. For each $b\in P_a$ such that $Q_b(c)\ne0$, the line through $c$ 
and a point $z\in C_b$ meets $C_b$ a second time, counting multiplicity.
The map $C_b\to C_b$ taking $z$ to the second intersection point is the reflection in $c$. 
In the plane $\{b\}\times P_z$, take $(b,c)$ as the center of reflection. Over the open subset of $P_a$ where $Q_b(c)\ne 0$, these reflections define a rational map $r_c:X\dasharrow X$ preserving the first projection $f:X\to P_a$. 
\end{notation}
\begin{lemma}\label{lemma:rc}
In homogeneous coordinates, $r_c$ is given by 
\begin{equation}\label{eq:rc}
 r_c(b,[z])=(b,[R_{b,c}(z)]),\quad\text{where}\quad
 R_{b,c}(z)=Q_b(c)z-2B_b(z,c)c.
\end{equation} 
It is a biregular involution on the open subset of $X$ where $Q_b(c)\ne 0$.
\end{lemma}
\begin{proof} We parametrize the line $\ell=zc\subset P_z$ by $[s:t]\mapsto[sz+tc]$. 
On $\ell$, the equation of $C_b$ becomes
\[
 t\bigl(2sB_b(z,c)+tQ_b(c)\bigr)=0.
\]
The root $t=0$ gives $z$. The residual root is $[s:t]=[Q_b(c):-2B_b(z,c)]$. 
Its image under the parametrization of $\ell$ is
\[
[sz+tc]=[Q_b(c)z-2B_b(z,c)c]=[R_{b,c}(z)].
\]
Since $R_{b,c}(R_{b,c}(z))=Q_b(c)^2z$, the map $r_c$ is a biregular involution 
on the open subset of $X$ where $Q_b(c)\ne 0$.
\end{proof}

Recall that $D_q\subset P_a$ stands for the fiber of the second conic bundle 
over a point $q\in P_z$, 
see \eqref{eq:conics}.

\begin{lemma}\label{chord}
Let $x=(b,p)$ be a smooth point of the conic $C_b$. 
If $C_b$ is smooth, set $M=C_b$. 
Otherwise, let $M$ be the component of $C_b$ that does not contain $p$. 
For $q$ in a nonempty open subset $\omega$ of $M$, 
the conic $D_q$ is smooth and every center $c\in pq\setminus C_b$ 
has the following properties: 
the reflection $r_c$ exchanges $x$ and $y=(b,q)$, is biregular near both points, 
and is regular 
along the whole of $D_q$.
\end{lemma}
\begin{proof}
We choose $q$ so that the line $\ell=pq$ is a component 
of no fiber of the conic bundle $f:X\to P_a$. 
If $C_a$ is singular and $p\in C_a$, then $a\in D_p\cap\Delta_f$. 
This set is finite, since $D_p$ has degree two and $\Delta_f$ is an irreducible sextic, 
see Lemma \ref{discriminants}. 
There are therefore only finitely many lines  through $p$ lying on fibers of the first conic bundle. 
Each such line meets $M$ in a finite set. If $M$ is a conic this follows from its irreducibility, 
and if $M$ is a line it misses $p$. Remove all these intersections from $M$. 
Also remove $M\cap\Delta_g$, which is finite by the same degree argument. 
If $C_b$ is smooth, remove $p$ as well. 
Only finitely many points have been removed. For each remaining $q$, 
the conic $D_q$ is smooth and $\ell$ cuts $C_b$ at the two distinct points $p,q$.

Take any $c\in\ell\setminus C_b$. Since $Q_b(c)\ne0$, the reflection 
exchanges $x$ and $y$ and is biregular near both. 
Now take any $(a,q)\in D_q$. If $R_{a,c}(q)$ vanished, 
the linear independence of representatives of $q$ and $c$ would give $Q_a(c)=B_a(q,c)=0$. 
We already have $Q_a(q)=0$. Thus
\[
 Q_a(sq+tc)=s^2Q_a(q)+2stB_a(q,c)+t^2Q_a(c)=0
 \qquad\text{for all }s,t.
\]
Thus $Q_a$ would vanish identically on $\ell$, making $\ell=pq$ a component of the fiber $C_a$. This contradicts our choice of $q$. 
Hence the formula \eqref{eq:rc} for $r_c$ defines a morphism near every point of $D_q$.
\end{proof}

\subsection{Three spray velocities}

\begin{definition} Let $s:U\times\mathbb A^1\to X$ be a rank-one spray 
on an open subset $U\subset X$ with values in $X$. 
We say that this spray is \emph{$f$-relative} if $f(s(x,t))=f(x)$ 
for every $(x,t)\in U\times\mathbb A^1$. 
The same definition applies to the second conic bundle $g$.
\end{definition}

\begin{lemma}\label{relative}
Let $x=(b,p)$ 
 be a smooth point of the conic $C_b$. There is an $f$-relative spray 
 on a neighborhood of $x$ in $X$ whose velocity at $x$ is nonzero. 
 The analogous statement applies to the second conic bundle $g:X\to P_z$. 
\end{lemma}
\begin{proof}
If $C_b$ is smooth, we choose $q\in C_b\setminus\{p\}$. If $C_b$ is reducible, 
we choose $q$ on the component missing $p$, away from the node. 
Then the line $\ell=pq$ meets $C_b$ transversely at $p$ and $q$. 
We fix $c\in\ell\setminus C_b$, so that $r_c(x)=(b,q)$. 
Since $c\notin C_b$, we have $Q_b(c)\ne0$. Thus the locus
$\{(a,z)\in X\mid Q_a(c)\ne0\}$ is open in $X$ and contains $x$.
We choose an open neighborhood $U$ of $x$ in this locus. For $u\in U$, we set $b_u=f(u)$ and write $u=(b_u,p_u)$ 
and $r_c(u)=(b_u,q_u)$. After shrinking $U$, we may assume that 
$p_u$ and $q_u$ are distinct and that the line $p_uq_u$ meets $C_{b_u}$ 
transversely at both points.

We write $P_z=\mathbb P(V_z)$, where $V_z$ is a three-dimensional
complex vector space. For $u=(b_u,p_u)\in U$, we identify the plane
$\{b_u\}\times P_z$ containing the fiber with $P_z$ by the second
projection. On this plane we choose homogeneous coordinates
$[Z_0:Z_1:Z_2]$, varying regularly with $u$, such that
\begin{equation}\label{eq:marked-fiber-coordinates}
 q_u=[1:0:0],\qquad p_u=[0:1:0],\qquad
 \mathbb T_{q_u}C_{b_u}=\{Z_1=0\}.
\end{equation}
Here $\mathbb T_{q_u}C_{b_u}$ denotes the projective tangent line to $C_{b_u}$
at $q_u$. To obtain these coordinates, we choose local vector
representatives $\widetilde q_u,\widetilde p_u\in V_z$ of $q_u,p_u$.
The projective tangent line is $\mathbb P(K_u)$, where
\[
 K_u=\ker\bigl(B_{b_u}(\widetilde q_u,-):V_z\to\mathbb C\bigr).
\]
Since $q_u$ is a smooth point of $C_{b_u}$, these kernels form a
rank-two subbundle of $U\times V_z$. After shrinking $U$, we complete
$\widetilde q_u$ to a local frame of this subbundle by a vector $v_u$.
Transversality gives $\widetilde p_u\notin K_u$, so the ordered basis
$(\widetilde q_u,\widetilde p_u,v_u)$ gives the coordinates
in~\eqref{eq:marked-fiber-coordinates}.

We now express the quadratic form $Q_{b_u}(z)=F(b_u,z)$ in these coordinates.
Since it vanishes at the two marked points
in~\eqref{eq:marked-fiber-coordinates}, its coefficients of $Z_0^2$
and $Z_1^2$ vanish. Since the tangent line to $C_{b_u}$ at $q_u=[1:0:0]$
has equation $Z_1=0$, the coefficient of $Z_0Z_2$ vanishes and the
coefficient of $Z_0Z_1$ is a unit on $U$. Dividing $Q_{b_u}$ 
by the coefficient of $Z_0Z_1$, we obtain the equation of the conic 
$C_{b_u}$ in the coordinates chosen at $u$:
\[
 Z_0Z_1+\beta(u)Z_1Z_2+\gamma(u)Z_2^2=0
\]
with $\beta,\gamma$ regular on $U$. In the chart $Z_1=1$, 
the conic is the graph $Z_0=-\beta(u)Z_2-\gamma(u)Z_2^2$. 
We define the spray by
\[
 s(u,t)=\bigl(b_u,[-\beta(u)t-\gamma(u)t^2:1:t]\bigr),
 \qquad (u,t)\in U\times\Aff^1.
\]
Then $s(u,0)=(b_u,p_u)=u$ and $f(s(u,t))=b_u=f(u)$. 
Its velocity is nonzero because $Z_2=t$ in this chart. 
If $C_{b_u}$ is reducible, then $\gamma(u)=0$ and the same graph 
parametrizes the affine part of the component through $p_u$.
\end{proof}

\begin{proof}[Proof of Theorem~\ref{thm:biquadratic}]
 
\smallskip\noindent\emph{Choosing the first spray.}
We fix $x=(b,p)\in X$ and choose affine coordinates centered at $b$ in $P_a$
and at $p$ in $P_z$. In the resulting product chart, we write $F=0$ for the equation of $X$.
For $\xi\in T_bP_a$ and $\eta\in T_pP_z$, we write
\[
d_aF_x(\xi)=dF_x(\xi,0),\qquad d_zF_x(\eta)=dF_x(0,\eta).
\]
Thus $dF_x(\xi,\eta)=d_aF_x(\xi)+d_zF_x(\eta)$. Since $X$ is smooth at $x$,
at least one of the two linear forms $d_aF_x$ and $d_zF_x$ is nonzero.
After interchanging the two planes if necessary, we may assume that $d_zF_x\ne0$.
Then $p$ is a smooth point of $C_b$. For every $\xi\in T_bP_a$, we can choose
$\eta\in T_pP_z$ such that 
\[dF_x(\xi,\eta)=d_aF_x(\xi)+d_zF_x(\eta)=0,\] 
since $d_zF_x\ne0$.
Then $(\xi,\eta)\in T_xX$ and $df_x(\xi,\eta)=\xi$, so $df_x$ is surjective.
Its kernel is the tangent line to $C_b$ at $p$. Lemma~\ref{relative} 
gives a first rank-one spray whose velocity at $x$ spans $\ker df_x$.
To obtain three velocities spanning $T_xX$, it remains to construct 
two more rank-one sprays near $x$ whose velocities at $x$ 
project under $df_x$ to a basis of $T_bP_a$.

\smallskip\noindent\emph{Choosing the two projected directions.}
Let $M$ and $\omega\subset M$ be as in Lemma~\ref{chord}. For each $q\in\omega$, 
the restriction of $f$ to the second fiber $D_q$ identifies it with the 
smooth plane conic $F(a,q)=0$. At $y=(b,q)$ its tangent projects isomorphically 
onto the line
\[
 T_bD_q=\ker\bigl(Q_u(q)\,du_b+Q_v(q)\,dv_b\bigr)\subset T_bP_a.
\]
The covector $Q_u(q)\,du_b+Q_v(q)\,dv_b\in T_b^*P_a$ is the differential at $b$
of the function $a\mapsto F(a,q)$. It is nonzero because $D_q$ is smooth at $b$.
If the tangent line $T_bD_q$ were the same for every $q\in\omega$, 
we could choose a fixed nonzero vector $\xi$ on it. 
Then $du_b(\xi)\,Q_u(q)+dv_b(\xi)\,Q_v(q)=0$
for every $q\in\omega$. The scalars $du_b(\xi)$ and $dv_b(\xi)$ cannot both
vanish, since $u,v$ are local coordinates and $\xi\ne0$.
The resulting linear combination of $Q_u,Q_v$ would vanish on all of $M$, 
contrary to Lemma~\ref{derivative-quadrics}. We may therefore choose $q_1,q_2\in\omega$ 
giving two distinct tangent lines. These lines span the two-dimensional space $T_bP_a$.

\smallskip\noindent\emph{Transporting a spray back to $x$.}
For either choice of $q\in\{q_1,q_2\}$, set $y=(b,q)$ and choose $r=r_c$ as in Lemma~\ref{chord}. 
Let $Z=X\setminus\Dom(r)$. The second projection $g:X\to P_z$ is projective and hence closed, so 
$V=P_z\setminus g(Z)$ is open. It contains $q$ because $r$ is regular on the entire fiber $D_q$. 
In particular $r$ is regular everywhere on $g^{-1}(V)$.

By Lemma~\ref{relative}, we choose a $g$-relative spray $s:U\times\mathbb A^1\to X$
on an open neighborhood $U$ of $y$, with nonzero velocity at $y$.
Since $y\in g^{-1}(V)$, we may shrink $U$ so that $U\subset g^{-1}(V)$.
For every $u\in U$ and $t\in\mathbb A^1$, we have $g(s(u,t))=g(u)\in V$.
Hence $s(u,t)\in g^{-1}(V)\subset\Dom(r)$.
We now shrink the neighborhood of $x$ so that $r$ is biregular there and maps it into $U$.
On this neighborhood
\[
 \widetilde s(w,t)=r\bigl(s(r(w),t)\bigr)
\]
is defined for every $t\in\Aff^1$ and satisfies $\widetilde s(w,0)=w$, 
because $r$ is an involution near $x$ and $y$.

Let $v$ be the velocity of $s$ at $y$. It spans $T_yD_q$. 
The velocity of $\widetilde s$ at $x$ is $dr_y(v)$, and differentiating $f\circ r=f$ at $y$ gives
\[
 df_x\bigl(dr_y(v)\bigr)=df_y(v).
\]
Thus the two transported velocities project onto the two independent directions chosen above. 
Together with the first velocity, which spans $\ker df_x$, they span $T_xX$. Since $x\in X$ was arbitrary, 
the local spray criterion of  \cite[Lemma~2.5]{KZcubic} gives algebraic ellipticity of $X$.

Lemmas~\ref{singular-fibers}--\ref{derivative-quadrics} supply 
a dense open set $\Omega\subset |\mathcal O(2,2)|$ consisting 
of smooth divisors $X\subset\PP^2\times\PP^2$ of bidegree  $(2,2)$
for which the construction works. By adjunction, we have $-K_X=\mathcal O_X(1,1)$. 
Hence these $X$  are Fano varieties.
By \cite[Theorem~1.1]{BB}, applied with $n=2$, a very general hypersurface of bidegree $(2,2)$ 
is stably irrational. Including the complement of $\Omega$ in the countable union 
of excluded closed subsets proves the last assertion.
\end{proof}

\section{Double covers of the incidence threefold}
Let $P=\PP^2$ and let $F\subset P\times P^\vee$ be the incidence threefold.
For $B\in H^0(F,\OO_F(2,2))$, let $X_B\to F$ be the double cover
branched along $S_B=\{B=0\}$, and let $\Delta_f,\Delta_g$ be the
discriminants of the projections to $P$ and $P^\vee$.

\begin{theorem}[Theorem~\ref{main:unified}(3)]\label{thm:cover}
\CoverStatement
\end{theorem}

\subsection{The two conic bundles}
Write $X=X_B$. In the total space of
$\OO_F(1,1)$, we define the double cover $\pi:X\to F$ by the equation
\[
 w^2=B(z),\qquad z\in F,\quad w\in\OO_F(1,1)_z.
\]
Here $w^2=w\otimes w$ belongs to the fiber of $\OO_F(2,2)$ at $z$.
The projection $\pi(z,w)=z$ is branched along $S_B$.
We denote the two projections by $f:X\to P$ and $g:X\to P^\vee$,
and their fibers by $C_b=f^{-1}(b)$ and $D_L=g^{-1}(L)$, respectively.
The discriminant schemes are defined locally by the determinants of the
fiberwise quadratic forms. Their underlying loci are
\[
 \Delta_f=\{b\in P\mid C_b\text{ is singular}\},\qquad
 \Delta_g=\{L\in P^\vee\mid D_L\text{ is singular}\}.
\]

We first verify the assertion about a general section $B$.

\begin{lemma}\label{fibers}
If $\Delta_f\ne P$ and $\Delta_g\ne P^\vee$, then $f$ and $g$ are
conic bundles. If the discriminant curves are smooth, every singular
fiber is the union of two distinct lines.
\end{lemma}
\begin{proof}
For $b\in P$, let $\Lambda_b=\{L\in P^\vee\mid b\in L\}$ be the pencil of
lines through $b$. The restriction of $\pi$ to $C_b$ is a double cover
$C_b\to\Lambda_b\simeq\PP^1$. We write $B_b$ for the restriction of the
branch section $B$ to $\{b\}\times\Lambda_b$. After choosing homogeneous
coordinates $[\lambda_0:\lambda_1]$ on $\Lambda_b$ and trivializing
$\OO_P(1)$ at $b$, we obtain the plane equation
\[
 C_b=\bigl\{[\lambda_0:\lambda_1:w]\in\PP^2
       \mid w^2=B_b(\lambda_0,\lambda_1)\bigr\}.
\]
Both projections are flat and projective, since each is the composite
of the finite flat double cover $\pi$ and a projective line bundle.
The displayed equation shows that $C_b$ is singular precisely when the binary quadratic form $B_b$
has zero discriminant. If $B_b$ has rank one, $C_b$ is the union of two
distinct lines. If $B_b=0$, it is a double line. Thus the general fiber
of $f$ is a smooth conic whenever $\Delta_f\ne P$.

If $B_b=0$, the discriminant, being quadratic in the coefficients of
$B_b$, vanishes to order at least two at $b$. Thus smoothness of
$\Delta_f$ excludes double lines. The same description applies to $g$
after interchanging the two planes.
\end{proof}

\begin{lemma}\label{general}
For a general section $B$, the surface $S_B$ is smooth and the two
discriminants $\Delta_f$ and $\Delta_g$ are smooth plane sextics.
\end{lemma}
\begin{proof}
Since $\OO_F(2,2)$ is globally generated, Bertini's theorem gives
smoothness of $S_B$ for general $B$.
Let $p:F\to P$ be the first projection. The branch section $B$ induces
a section of the rank-three vector bundle $\mathcal E=p_*\OO_F(2,2)$
whose value at $b$ is $B_b$. We fix $b\in P$ and choose a section $s$ of
$\OO_P(2)$ with $s(b)\ne 0$. Since $\Lambda_b$ is a line in $P^\vee$,
every section of $\OO_{\Lambda_b}(2)$ extends to a section $t$ of
$\OO_{P^\vee}(2)$. As $t$ varies, the sections $(s\boxtimes t)|_F$
therefore generate the fiber
\[
 \mathcal E_b=\OO_P(2)_b\otimes H^0(\Lambda_b,\OO_{\Lambda_b}(2)).
\]
We set $V=H^0(F,\OO_F(2,2))=H^0(P,\mathcal E)$. For every $b\in P$,
evaluation $V\to\mathcal E_b$ is surjective by the preceding argument.
Thus the pairs $(b,B)$ with $B_b=0$ form a vector bundle over $P$ of
rank $\dim V-3$, namely the kernel of the evaluation map of bundles.
Since $\dim P=2$, its total space has dimension $\dim V-1$. Its image
in $V$ is closed because $P$ is projective. A general $B\in V$
therefore satisfies $B_b\ne0$ for every $b\in P$.

We now consider the universal evaluation map
\[
 e:P\times V\longrightarrow\operatorname{Tot}(\mathcal E),
 \qquad (b,B)\longmapsto B_b.
\]
This map is smooth, since it is a surjective morphism of vector bundles
over $P$. The rank-one forms constitute a smooth hypersurface $R_1$
in the complement of the zero section. Indeed, in a local frame the
determinant is $ac-d^2$ for a form with coefficients $(a,2d,c)$,
and its differential vanishes only at the zero form. Hence $e^{-1}(R_1)$
is smooth. By generic smoothness, its projection to $V$ is smooth over
a dense open subset of $V$. For a general $B$ as above, the fiber of
this projection is precisely the discriminant $\Delta_f$, with its
local determinant equation. Thus the section meets $R_1$ transversely
and $\Delta_f$ is smooth.

The fiber of
$p|_{S_B}:S_B\to P$ over $b$ is the zero scheme of $B_b$ on
$\Lambda_b$. Since $B_b$ is nonzero for every $b$, this map is a
double cover with branch curve $\Delta_f$. Adjunction gives
$K_F=\OO_F(-2,-2)$ and hence $K_{S_B}=\OO_{S_B}$. The canonical
bundle formula for the double cover then gives $\deg\Delta_f=6$.
Interchanging the two planes gives the same conclusion for $\Delta_g$.
Both conclusions hold on the intersection of the two dense open sets
of branch sections.
\end{proof}

By Lemma~\ref{general}, the hypotheses of Theorem~\ref{thm:cover} hold for
a general section $B$. Stable irrationality for very general $B$ follows
from \cite[\S7.2, following Lemma~7]{HTF}. It remains to prove
ellipticity under the stated smoothness hypotheses.
By the local criterion recalled in the introduction, it suffices to span
each tangent space by velocities of local rank-one sprays.

\subsection{Fiberwise reflections at a smooth point of a conic}
\begin{proposition}\label{prop:involution}
Under the hypotheses of Theorem~\ref{thm:cover}, let $x\in X$ be a smooth
point of at least one of the two conics $C_{f(x)}$ and $D_{g(x)}$.
Then there exist three local rank-one sprays at $x$ such that their velocities at $x$ span $T_xX$.
\end{proposition}

By symmetry we may assume below that $x$ is smooth on $C_b$, where $b=f(x)$.

\begin{lemma}\label{lem:relative}
There exists a local $f$-relative spray of rank one at $x$ whose
velocity spans $T_xC_b=\ker df_x$.
\end{lemma}
\begin{proof}
We trivialize the ambient plane bundle near $b$, writing its planes as
$\PP(V)$, where $\dim V=3$. We choose a chord through $x$ meeting $C_b$
transversely at a second smooth point $y$. If $C_b$ is reducible, we
choose $y$ on the other component. We fix a point $e\in\overline{xy}\setminus C_b$.
In the chosen trivialization over a neighborhood $U$ of $b$, we extend
$e$ to a section $\varepsilon$ of the plane bundle by keeping its
homogeneous coordinates fixed. Thus $\varepsilon(b)=e$, and
$\varepsilon(a)$ lies in the plane of $C_a$ for each $a\in U$. After
shrinking $U$, we have $\varepsilon(a)\notin C_a$ for every $a\in U$.
Reflection with center $\varepsilon(a)$ in the plane of $C_a$ defines
a regular involution of $f^{-1}(U)$, hence a birational involution
$r_\varepsilon\colon X\dashrightarrow X$. On $C_b$ it exchanges $x$
and $y$. Thus, for $u$ near $x$,
the point $v=r_\varepsilon(u)$ varies regularly with $u$ and lies on the same
conic $C_{f(u)}$.

In the chosen trivialization the conics form a family of plane conics,
and the chord from $u$ to $v$ remains transverse at both points after
shrinking about $x$. These are precisely the data used in the proof of
Lemma~\ref{relative}: take $u$ and $v$ as the two marked points in
\eqref{eq:marked-fiber-coordinates}. The same affine parametrization
therefore gives an $f$-relative spray with nonzero velocity at $x$,
including when $C_b$ is reducible. By symmetry the construction applies to $g$.
\end{proof}

\begin{lemma}\label{lem:involution}
For a general line $L$ through $b$, there exist a point $y\in C_b\cap D_L$
and a fiberwise reflection $r$ of $X$ over $P$ with the following
properties. The conic $D_L$ is smooth, and $df_y$ maps $T_yD_L$
isomorphically onto $T_bL$. The fiberwise reflection $r$ exchanges $x$ and $y$, is
biregular near both points, and is regular along the whole conic $D_L$.
\end{lemma}
\begin{proof}
We choose a general line $L$ through $b$, different from the line coordinate
of $x$, and a point $y\in C_b$ above $L$. In the reducible case choose
the component missing $x$; each component maps isomorphically to the
pencil. We may require $L\notin\Delta_g$ and $B(b,L)\ne0$.
Indeed, the pencil is a line in the dual plane and $\Delta_g$ is a
smooth sextic. Thus $D_L$ is smooth, and $df_y$ identifies its tangent
with $T_bL$ because the double cover is unramified at $y$.

We construct a fiberwise reflection $r$ over $P$ that exchanges $x,y$ and is
regular along $D_L$. In the plane of $C_b$ choose a center $e$ on the
chord $xy$ outside $C_b$ whose projection to the pencil is a line
$m\ne L$. This is possible since $x,y$ project to distinct lines and
their chord is not a component. Choose $c\in m\setminus L$. For $a\ne c$
set $m(a)=ac$. This gives a section of the incidence bundle, and
$m(a)\ne L$ for all $a\in L$. Its graph pulls $\OO_F(1,1)$ back to $\OO_P(2)$
on $P\setminus\{c\}$, since $m(a)$ has linear coordinates in $a$.
To choose a point $e(a)$ in the plane of $C_a$ projecting to $m(a)$,
we specify a vector $w\in\OO_F(1,1)_{(a,m(a))}$. Under the preceding
identification, a section $h\in H^0(P,\OO_P(2))$ gives a moving center
$e(a)$ by setting $w=h(a)$. We prescribe $h(b)$ so that $e(b)=e$.
For $a\in(L\cap\Delta_f)\setminus\{b\}$ require $e(a)\notin C_a$.
At each of these finitely many points this excludes at most two values,
given by $h(a)^2=B(a,m(a))$ in a trivialization. For each such $a$,
we choose a linear form $\ell$ vanishing at $b$ and nonzero at $a$.
Adding a scalar multiple of $\ell^2$ to $h$ preserves $h(b)$ and allows
$h(a)$ to take any value in $\OO_P(2)_a$. Thus each forbidden value
cuts out a proper affine hyperplane in the space of quadratic forms
with the prescribed value at $b$. We choose $h$ outside the union of
these finitely many hyperplanes.

Let $Q_a$ denote the ternary equation of $C_a$ and $\beta_a$ its polar form.
In local frames the reflection is given by
\[
 z\longmapsto Q_a(e(a))z-2\beta_a(z,e(a))e(a).
\]
This homogeneous expression is invariant projectively under changes
of frame and scaling of the equation and of the center. It exchanges
the ends of a chord and its square on vectors is multiplication by
$Q_a(e(a))^2$. In particular it is biregular near $x$ and $y$.
For $z\in D_L$ over $a$ the projective points $z,e(a)$ are distinct,
as their line coordinates are $L,m(a)$. If the displayed vector vanished,
then $Q_a(e(a))=\beta_a(z,e(a))=0$, so the entire line joining $z,e(a)$
would lie in $C_a$. This is impossible for a smooth conic, and on
each singular conic over $L$ we arranged $e(a)\notin C_a$.

\end{proof}

\begin{proof}[Proof of Proposition~\ref{prop:involution}]
We choose $L,y,r$ as in Lemma~\ref{lem:involution}. Since $g$ is proper,
$V=P^\vee\setminus g(X\setminus\operatorname{Dom}(r))$ is an open
neighborhood of $L$. By Lemma~\ref{lem:relative} applied to $g$ at $y$,
there is a $g$-relative spray $s:W\times\Aff^1\to X$ whose velocity
spans $T_yD_L$. We choose a neighborhood $U$ of $x$ on which $r$ is
regular and such that $r(U)\subset W\cap g^{-1}(V)$. For every $u\in U$
and $t\in\Aff^1$, the point $s(r(u),t)$ lies on $D_{g(r(u))}$ by
$g$-relativity. Since $r$ is regular on $g^{-1}(V)$, the formula
\[
 \widetilde s:U\times\Aff^1\longrightarrow X,\qquad
 \widetilde s(u,t)=r\bigl(s(r(u),t)\bigr)
\]
defines a morphism on all of $U\times\Aff^1$, with
$\widetilde s(u,0)=r(r(u))=u$. Since $f\circ r=f$, the velocity of
$\widetilde s$ at $x$ projects to a generator of $T_bL$. Two distinct
general lines through $b$ thus give two velocities whose projections
span $T_bP$. Together with the velocity spanning $\ker df_x$ from
Lemma~\ref{lem:relative}, they span $T_xX$.
\end{proof}

\subsection{The associated del Pezzo fibration}\label{subsec:pencil}
We study a point $x\in X$ that is a node of both conics through it.
We choose a general line $\ell$ through $b=f(x)$ and a second point
$d\in\ell\cap\Delta_f$ distinct from $b$. The pencil through $d$
identifies a neighborhood of $x$ with an open subset of a family of
surfaces $\mathcal S\to J$. In
Proposition~\ref{prop:pencil} we realize $\mathcal S$ as two successive
blowups of a fibration $\mathcal T\to J$ in del Pezzo surfaces of degree
four, first along a section and then along a section of the resulting
exceptional divisor. The contraction $\mathcal S\to\mathcal T$ is an
isomorphism near $x$.

We retain the hypotheses of Theorem~\ref{thm:cover} and write $\pi(x)=(b,L)$.
For a line $m\subset P$, we denote by $m^\perp$ the corresponding
point of $P^\vee$. We choose a general line $\ell\subset P$ through $b$,
distinct from $L$, meeting $\Delta_f$ transversally and satisfying
$\ell^\perp\notin\Delta_g$.

\begin{lemma}\label{lem:pencil-tangent}
The maps $df_x$ and $dg_x$ have rank one, 
and $\operatorname{im}df_x=T_bL$.
For a general line $\ell$ through $b$, the surface $S_\ell=f^{-1}(\ell)$
is smooth and $T_xS_\ell=\ker df_x$.
\end{lemma}
\begin{proof}
Let $y=\pi(x)=(b,L)\in F$. Locally on $F$ near $y$, the double cover 
$\pi\colon X\to F$ has equation
$w^2=B$, where $B$ is a regular function whose zero locus is $S_B$. 
The two ruling lines through
$y$ are obtained by fixing $b$ and varying the line through it, or by
fixing $L$ and varying the point on it. For either ruling line $R$,
choose a local parameter $t$ on $R$ vanishing at $y$. The inverse
image of $R$  in $X$ has equation $w^2=B|_R(t)$ 
and is singular at $x$ by
hypothesis. The two partial derivatives therefore vanish at $x$,
giving $w(x)=0$ and $(B|_R)'(0)=0$. The latter equality says that
$dB_y$ annihilates $T_yR$. Thus both ruling directions belong to
$T_yS_B=\ker dB_y$. They are independent, so they span this plane.
Differentiating $w^2=B$ at $w=0$ gives
$\operatorname{im}d\pi_x=\ker dB_y$. Under the projection to $P$,
the first ruling line maps to $b$ and the second maps isomorphically
onto $L$. Consequently
$\operatorname{im}df_x=T_bL$, and both 
$df_x$ and $dg_x$ have rank one.

The pencil of lines through $b$ has no base point on
$\Delta_f\setminus\{b\}$. By Bertini's theorem, a general member $\ell$
of this pencil meets $\Delta_f$ transversely away from $b$. Since
$\Delta_f$ is smooth at $b$, transversality at $b$ amounts to choosing
$\ell\ne\mathbb T_b\Delta_f$. Moreover, the line $\Lambda_b\subset P^\vee$
parametrizing the pencil meets the smooth sextic $\Delta_g$ in at most
six points. We may therefore choose a line $\ell$ through $b$ 
meeting $\Delta_f$ transversally
at every point of intersection, with $\ell\ne L$ and $\ell^\perp\notin\Delta_g$.
For each of the six points $a\in\ell\cap\Delta_f$, we denote by $y_a$
the unique node of the singular conic $C_a=f^{-1}(a)$. We choose local
analytic coordinates $(u,v,t)$ on $X$ at $y_a$ and $(s,t)$ on $P$ at
$a$ in which $f$ is given by $f(u,v,t)=(uv,t)$ and $\Delta_f$ has
equation $s=0$. Since $\ell$ meets $\Delta_f$ transversally at $a$,
we may write $\ell$ locally as $t=h(s)$. Its inverse image $S_\ell$
then has equation $t=h(uv)$, so it is smooth at $y_a$.

At every point $z\in S_\ell$ distinct from these six nodes, the conic
$C_{f(z)}$ is smooth at $z$. Hence $f$ is smooth at $z$, and so is
$S_\ell=f^{-1}(\ell)$. Since the distinct
lines $\ell$ and $L$ have distinct tangent directions at $b$, we obtain
$T_xS_\ell=(df_x)^{-1}(T_b\ell)=\ker df_x$.
\end{proof}

Since $\Delta_f$ is a plane sextic,
$\ell\cap\Delta_f$ consists of six distinct points. One of them is
$b$. We fix one of the other five and denote it by $d$.
We write $A$ for either component of $C_d$ and set $R=D_{\ell^\perp}$.

\begin{lemma}\label{lem:pencil-contractions}
The curves $A$ and $R$ on $S_\ell$ are smooth rational curves with
$A^2=-1$, $R^2=-2$, and $A\cdot R=1$. Contracting $A$ and then
the image of $R$ gives a smooth del Pezzo surface of degree four.
\end{lemma}
\begin{proof}
On $S_\ell$ the two components of the conic $C_d$ meet transversely once, so
$A^2=-1$. Restricting the incidence variety $F$ over the line $\ell\subset P$,
we obtain the surface
\[
F_\ell=F\cap(\ell\times P^\vee)
=\{(p,m^\perp)\in\ell\times P^\vee\mid p\in m\}.
\]
The double cover $\pi\colon X\to F$ restricts to a double cover
$\pi_\ell\colon S_\ell\to F_\ell$. For $m\ne\ell$, the unique point of
$F_\ell$ over $m^\perp$ is $(m\cap\ell,m^\perp)$. Thus $F_\ell$ is the
closure of the graph of the linear projection $P^\vee\dashrightarrow\ell$
given by $m^\perp\mapsto m\cap\ell$. Its second projection is therefore
the blowup of $P^\vee$ at $\ell^\perp$, so $F_\ell$ is the Hirzebruch
surface $\mathbb F_1$.

The exceptional curve $E=\ell\times\{\ell^\perp\}$ is the negative
section of the first projection $F_\ell\to\ell$. Since $g$ records the
line in each incident pair, $\pi_\ell^*E=g^{-1}(\ell^\perp)=R$.
This is a smooth conic because $\ell^\perp\notin\Delta_g$, and the
degree-two intersection formula gives $R^2=2E^2=-2$. Moreover $A\cdot R=1$, since the map from $A$ to the
pencil through $d$ is an isomorphism and $R$ selects the member $\ell$.
After contraction of $A$, the image of $R$ is smooth and has square
$-2+1=-1$. We contract this image as well and denote the resulting smooth
projective surface by $T$. Writing $H_1,H_2$ for the pullbacks to $X$
of the hyperplane classes of $P,P^\vee$, respectively, we obtain by
adjunction $-K_{S_\ell}=H_2|_{S_\ell}$ and $K_{S_\ell}^2=H_1H_2^2=2$. Since $R$ is the
only curve contracted by $g|_{S_\ell}$, every other irreducible curve on $S_\ell$
has positive degree with respect to $-K_{S_\ell}$.

Let $C\subset T$ be an irreducible curve. Its strict transform on $S_\ell$
is distinct from $R$, hence has positive degree with respect to $-K_{S_\ell}$.
At each step, let $E$ be the $(-1)$-curve we contract and let $D$ be
the strict transform of $C$ on the same surface. We obtain the
anticanonical degree of the image of $D$ by adding $D\cdot E$ to the
anticanonical degree of $D$. Since $D\cdot E\geq 0$, the degree stays
positive through both contractions, so $(-K_T)\cdot C>0$.
Moreover, each contraction increases the square of the canonical class
by one, giving $K_T^2=K_{S_\ell}^2+2=4$. By the Nakai--Moishezon criterion,
$-K_T$ is ample. Thus $T$ is a smooth del Pezzo surface of degree four.
\end{proof}

\begin{proposition}\label{prop:pencil}
With the notation of this section, the following assertions hold.
\begin{enumerate}
\item The surface $S_\ell=f^{-1}(\ell)$ is smooth, and $T_xS_\ell=\ker df_x$.
\item There is an open neighborhood $J$ of $\ell^\perp$ in the pencil
of lines through $d$ for which the following holds. Let $\mathcal S\to J$
be the incidence family of surfaces $S_m=f^{-1}(m)$. Then there exist a smooth
projective family $\mathcal T\to J$ of del Pezzo surfaces of degree
four, a morphism $\varphi:\mathcal S\to\mathcal T$ over $J$, and
a section $\sigma:J\to\mathcal T$ such that $\varphi$ is an isomorphism
over $\mathcal T\setminus\sigma(J)$. Moreover,
$\varphi(x,\ell^\perp)\notin\sigma(J)$.
\end{enumerate}
\end{proposition}
\begin{proof}[Proof of Proposition~\ref{prop:pencil}]
Lemma~\ref{lem:pencil-tangent} proves assertion~(1).
To prove (2), consider the pencil of lines through $d$, parametrized by
the dual line $d^\perp=\Lambda_d\subset P^\vee$. Let $J\subset d^\perp$
be the open subset parametrizing the lines $m$ through $d$ 
that meet $\Delta_f$
transversally and satisfy $m^\perp\notin\Delta_g$. By our choice of
$\ell$, we have $\ell^\perp\in J$. Let $\beta\colon\widetilde P\to P$ 
be the blowup at $d$.
The pencil of lines through $d$ induces the ruling
$q\colon\widetilde P\to\Lambda_d$ of the Hirzebruch surface
$\widetilde P\simeq\mathbb F_1$. Pulling back $f\colon X\to P$ by
$\beta$ and composing with $q$, we obtain a family of surfaces over
$\Lambda_d$. Its fiber over $m^\perp$ is $S_m=f^{-1}(m)$, since
$\beta$ maps the corresponding ruling fiber isomorphically onto $m$.
Over $J$, this gives the incidence family
\[
 \mathcal S=\{(u,m^\perp)\in X\times J\mid f(u)\in m\}
 \longrightarrow J,\qquad (u,m^\perp)\longmapsto m^\perp.
\]
Its fiber over $\ell^\perp$ is $S_\ell$.
The family $\mathcal S\to J$ is projective, since $\mathcal S$ is closed in
$X\times J$. It is flat, since the incidence equation cuts out a Cartier
divisor on every fiber of $X\times J\to J$. Each $m$ parametrized by $J$
meets $\Delta_f$ transversely, so the argument of Lemma~\ref{lem:pencil-tangent}
shows that $S_m$ is smooth. Thus $\mathcal S\to J$ is smooth.

Consider the divisors $\mathcal A=A\times J$ and
$\mathcal R=\{(u,m^\perp)\in\mathcal S\mid g(u)=m^\perp\}$.
For each $m^\perp\in J$, the fiber of $\mathcal R\to J$ is the smooth
conic $D_{m^\perp}=g^{-1}(m^\perp)$. The restriction of $f$ gives a double cover
\[
 f|_{D_{m^\perp}}\colon D_{m^\perp}\longrightarrow m.
\]
Thus $H_1$ has degree two on $D_{m^\perp}$, while $H_2$ restricts
trivially because $g$ is constant there. The polarization $H_1+H_2$
therefore restricts to $\mathcal O_{\mathbb P^1}(2)$ on this fiber.

On each surface $S_m$, the curves $A$ and $D_{m^\perp}$ meet transversely
at the unique point of $A$ lying over $m^\perp$ under $g$. Since
$g|_A:A\to\Lambda_d$ is an isomorphism, these points form a section
over $J$. Thus $\mathcal R\to J$ is a smooth family of conics with a
section, and both $\mathcal A$ and $\mathcal R$ are ruled surfaces over $J$.

We next construct $\mathcal T$ by two successive divisorial Mori contractions over $J$.
Contracting $A$ in each fiber gives a morphism
$\theta:\mathcal S\to\mathcal S'$ to a smooth family of surfaces over $J$.
It maps $\mathcal A$ onto a section of this family. Since $D_{m^\perp}$
meets $A$ transversely once, $\theta$ maps $D_{m^\perp}$ isomorphically
onto a smooth rational curve of self-intersection $-2+1=-1$.
Thus $\theta$ identifies $\mathcal R$ with a ruled surface in $\mathcal S'$.
For each $m^\perp\in J$, its ruling fiber $\theta(D_{m^\perp})$ has
self-intersection $0$ on $\theta(\mathcal R)$ and self-intersection $-1$
on the corresponding surface $\mathcal S'_{m^\perp}$. Contracting these curves in each fiber of 
$\mathcal S'\to J$ gives a second
morphism $\tau:\mathcal S'\to\mathcal T$, which maps $\theta(\mathcal R)$
onto the image of a section $\sigma:J\to\mathcal T$. We have
\[
 \mathcal S\xrightarrow{\theta}\mathcal S'\xrightarrow{\tau}\mathcal T,
 \qquad \varphi=\tau\circ\theta.
\]
Here $\mathcal T\to J$ is a smooth family of del Pezzo surfaces of degree
four, and $\varphi$ is an isomorphism over $\mathcal T\setminus\sigma(J)$.

The existence of both contractions over $J$ follows from \cite[Lemma~2.6(a)]{FKZfamilies},
applied first to $\mathcal A$ in the relative normal crossings divisor
$\mathcal A+\mathcal R$, and then to $\theta(\mathcal R)$ in $\mathcal S'$.
The resulting families are proper and flat over $J$, and their fibers are
smooth by construction. By Lemma~\ref{lem:pencil-contractions}, each fiber
of $\mathcal T\to J$ is a del Pezzo surface of degree four. Hence the
relative anticanonical bundle is relatively ample, so $\mathcal T\to J$
is projective. Finally, $x$ lies on neither $A$ nor $R$, since $b\ne d$
and $L\ne\ell$. Thus $\varphi(x,\ell^\perp)\notin\sigma(J)$.
\end{proof}

\subsection{Sprays avoiding the distinguished section}
By Proposition~\ref{prop:involution}, if $x\in X$ is smooth on at least one
of the two conics $C_{f(x)}$ and $D_{g(x)}$ through it, there are three
local rank-one sprays with values in $X$ whose velocities at $x$ span
$T_xX$. It remains to construct a dominating family of local sprays near
a point $x$ singular on both conics. For this we use
Proposition~\ref{prop:pencil} and retain the notation of Subsection~\ref{subsec:pencil}.

Let $a=\varphi(x,\ell^\perp)\in\mathcal T$. We shall construct two local
sprays near $a$ that preserve the fibers of $\mathcal T\to J$ and whose
images avoid $\sigma(J)$. Since $\varphi$ is an isomorphism over
$\mathcal T\setminus\sigma(J)$, each spray lifts to $\mathcal S$.
Composing with the projection $\mathcal S\to X$, which is an isomorphism
near $(x,\ell^\perp)$, gives two local sprays with values in $X$. Their
velocities at $x$ span $\ker df_x$. Interchanging $f$ and $g$ gives two
more sprays whose velocities span $\ker dg_x$. Together these velocities
span $T_xX$.

Set $j_0=\ell^\perp\in J$ and $T_0=\mathcal T_{j_0}$. On the fiber over
$j_0$, the morphism $\varphi$ sends $x\in S_\ell$ to the point
$a=\varphi(x,j_0)\in T_0$. Let $z=\sigma(j_0)$ be the point where the
distinguished section meets $T_0$. By Proposition~\ref{prop:pencil}, $a\ne z$.

\begin{lemma}\label{lem:nodal-avoidance}
For any tangent direction $\delta$ at $a$ in $T_0$,
there exist a point $q\ne a$ on no line of $T_0$ 
in its anticanonical embedding
and an integral curve $C\in|-K_{T_0}|$ through $a$ with the following
properties. Its only singularity is an ordinary node at $q$, it avoids $z$,
and $T_aC$ has direction different from $\delta$.

\end{lemma}
\begin{proof}
Apply \cite[Lemma~2.1]{DKZ26b} to $T_0$ with forbidden set $\{z\}$.

\end{proof}

\begin{lemma}\label{lem:relative-conic-pencil}
Let $q$ and $C$ be as in Lemma~\ref{lem:nodal-avoidance}. 
After shrinking $J$ about $j_0$,
there is a section $\gamma:J\to\mathcal T$ through $q$. The blowup
$\mu:\mathcal M\to\mathcal T$ along $\gamma(J)$ admits a proper flat
morphism $k:\mathcal M\to\mathcal P$ over $J$, where $\mathcal P\to J$
is a projective line bundle. For $j\in J$, the fibers of $k$ over
$\mathcal P_j$ are the conics obtained from the anticanonical pencil on
$\mathcal T_j$ with a double point at $\gamma(j)$. In particular, the
strict transform of $C$ is a smooth fiber of $k$.

\end{lemma}
\begin{proof}
The generic fiber of $\mathcal T\to J$ is rationally connected
and the family has good reduction at $j_0$. Weak approximation at $j_0$ \cite{HT}
gives the section through $q$, as in \cite[Lemma~3.2]{DKZ26b}.
We apply the construction in the proof of \cite[Proposition~3.6]{DKZ26b} over a
neighborhood of $j_0$. Blowing up the section gives a family of cubic
del Pezzo surfaces with a distinguished line in each fiber. The planes
through that line cut out the residual conics and form the required pencil.

\end{proof}

\begin{proof}[Proof of Theorem~\ref{thm:cover}]
Let $x\in X$ be a point singular on both conics through it.
We shall construct, on a neighborhood of $x$, a family of rank-one
sprays with values in $X$ whose velocities at $x$ span $T_xX$.
We choose $q$ and $C$ as in Lemma~\ref{lem:nodal-avoidance}, 
and we form the pencil $k$ of
Lemma~\ref{lem:relative-conic-pencil}. Let $\widetilde C$ 
be the strict transform of $C$ under $\mu$,
and let $p_0\in\mathcal P_{j_0}$ be the unique point such that
$k^{-1}(p_0)=\widetilde C$.
Let $V\subset\mathcal P$ be the open set over which the fibers of $k$ are
smooth, with the closed set $k(\mu^{-1}(\sigma(J)))$ removed. The latter set
is closed because $k$ is proper. Since $C$ avoids $z$, $p_0\in V$.
Thus $\mu(k^{-1}(p))$ avoids $\sigma(J)$ for every $p\in V$.

Choose an open neighborhood $U$ of $a$ disjoint from $\gamma(J)$.
Every $u\in U$ has a unique lift $\widetilde u$ under $\mu$. Shrinking $U$,
we obtain a morphism $h:U\to V$, $u\mapsto k(\widetilde u)$.
In the family $U\times_V k^{-1}(V)\to U$, the conic over $u$ is marked
at $\widetilde u$. These marked points form the section
$u\mapsto(u,\widetilde u)$. After shrinking $U$ about $a$, we identify
this family with $U\times\mathbb P^1$ so that the marked section becomes
the zero section $u\mapsto(u,0)$.
On the chart $U\times\mathbb A^1$, projection to $\mathcal M$ followed by
$\mu$ gives the desired morphism $s$ of the spray. 
For each $u$, the entire curve
$s(u,\mathbb A^1)$ avoids $\sigma(J)$ because it lies in
$\mu(k^{-1}(h(u)))$. Since $\mu$ is an isomorphism near $\widetilde a$,
the velocity of $s$ at $(a,0)$ is nonzero and spans $T_aC$.

We repeat the construction taking $\delta=T_aC$ 
in Lemma~\ref{lem:nodal-avoidance}, and we restrict the
two sprays to a common neighborhood of $a$. 
Their velocities span $T_aT_0$.

Both sprays can be lifted under $\varphi^{-1}$. 
The projection $\mathcal S\to X$ is an isomorphism
near $(x,\ell^\perp)$. Indeed, for $u$ near $x$, let $m(u)$ be the line
joining $d$ to $f(u)$. This line is well defined because $f(x)=b\ne d$.
Since $m(x)=\ell$ and $\ell^\perp\in J$, we have $m(u)^\perp\in J$
after shrinking the neighborhood of $x$. The inverse is
$u\mapsto(u,m(u)^\perp)$. Composing the lifted sprays 
with this projection gives two sprays on $X$ with velocities spanning
$\ker df_x$. We now interchange $f$ and $g$. Choose a general line
$n\subset P^\vee$ through $g(x)$ and a point $e\in n\cap\Delta_g$
distinct from $g(x)$. Applying the same construction to the pencil of
lines in $P^\vee$ through $e$, we obtain two further local rank-one
sprays with values in $X$. Their velocities at $x$ span
$T_x(g^{-1}(n))$, which equals $\ker dg_x$ by Lemma~\ref{lem:pencil-tangent}
with the two projections interchanged.
Since $(f,g)$ is $\pi$ followed by the incidence embedding,
$\ker df_x\cap\ker dg_x=\ker d\pi_x$ has dimension one. The two
planes therefore span $T_xX$. 
This completes the proof of
Theorem~\ref{thm:cover}.
\end{proof}

\section{Two conic bundles on a hypersurface in a product of rational surfaces}
\label{sec:rational-surfaces}

\begin{theorem}[Theorem~\ref{main:unified}(1)]\label{thm:rational-surfaces}
Let $S_1$ and $S_2$ be smooth rational projective surfaces over $\mathbb C$,
and let $X\subset S_1\times S_2$ be a smooth hypersurface. Suppose each projection
$X\to S_i$ is a conic bundle, with every geometric fiber isomorphic as a scheme
to a plane conic. Suppose that at least one surface is different from $\PP^2$,
or that both surfaces are planes and the bidegree of $X$ is either $(1,n)$ or $(n,1)$.
Then one of the projections admits a rational section, and $X$ is rational and algebraically elliptic.
\end{theorem}

\begin{proof}
Suppose first that both surfaces are planes. A general fiber of either projection is a smooth rational plane curve, so by the plane genus formula its degree is one or two. Thus $n\in\{1,2\}$ under our hypotheses. For bidegree $(1,n)$ the generic fiber of the second projection is a line defined over $\CC(S_2)$ and therefore has a rational point over that field. This gives a rational section and makes $X$ birational to $S_2\times\PP^1$. For bidegree $(n,1)$ we interchange the two projections.

Suppose now that at least one surface is different from $\PP^2$.
Let $C\subset S=S_2$ be a general fiber of the projection to $S_1$, viewed as a curve on $S_2$. It is a smooth rational curve. Its class $[C]\in{\rm Pic}(S)$ is nef: the curves move in a family covering $S$, so a general member meets any given irreducible curve properly. The intersection numbers $C\cdot H$, for $H\in\operatorname{Pic}(S)$, form a subgroup $\Lambda\subset\mathbb Z$. By adjunction, $\Lambda$ contains
\[
 C\cdot(C+K_S)=-2.
\]
Thus the positive generator $\tau$ of $\Lambda$ is either $1$ or $2$. If $\tau=1$, choose $H$ with $C\cdot H=1$. Pulling $H$ back from $S=S_2$ and restricting to the generic conic over $\mathbb C(S_1)$ gives a line bundle of degree one. The Riemann--Roch theorem supplies a nonzero section, whose zero divisor has degree one and hence gives a rational point. This is a rational section of the first projection. The same argument applies to the second projection.

Suppose that $\tau=2$. On a smooth rational surface the Picard group is a free unimodular intersection lattice (this can be seen by using a minimal rational model and the blowup formula). By unimodularity, there is a divisor class $D$ such that $D\cdot H=(C\cdot H)/2$ for every $H\in\operatorname{Pic}(S)$. By nondegeneracy, $C\sim 2D$, and $D$ is nef because $C$ is nef. We claim that the linear system $|D|$ defines a birational morphism $\varphi_{|D|}:S\to\mathbb P^2$ such that $D$ is the pullback of a line.

Set $d=D^2$. Adjunction for $C$ gives $K_S\cdot D=-2d-1$. Since $D^2+K_S\cdot D$ is even, $d$ is odd. Nefness gives $d\ge0$, so $d\ge1$ and $D$ is big and nef. Since $(K_S-D)\cdot D=-3d-1<0$, the class $K_S-D$ has no effective representative. Riemann--Roch and Serre duality give
\[
 h^0(S,\mathcal O_S(D))\ \ge\ \chi(\mathcal O_S(D))=(3d+3)/2.
\]
Since $D$ is nef and big, its $\Delta$-genus is nonnegative by \cite[Theorem~1.1]{Fujita}. The preceding Riemann--Roch estimate gives
\[
 0\leq\Delta(S,D)=2+d-h^0(S,\mathcal O_S(D))
 \leq 2+d-\frac{3d+3}{2}=\frac{1-d}{2}.
\]
Since $d\geq1$, we obtain $d=1$, $\Delta(S,D)=0$, and $h^0(S,\mathcal O_S(D))=3$. By the equality case of the same theorem, there exist a birational morphism $\pi:S\to W$ and a very ample line bundle $H$ on $W$ such that $\mathcal O_S(D)\simeq\pi^*H$. Since $H^2=D^2=1$, the embedding defined by $|H|$ realizes $W$ as a plane. Its three independent linear forms pull back to a basis of $H^0(S,\mathcal O_S(D))$. Thus $\pi:S\to\mathbb P^2$ is the morphism defined by $|D|$.

Suppose now that on each surface the class of a general fiber of the opposite projection is divisible by two. The claim gives birational morphisms $\pi_i:S_i\to\mathbb P^2$ such that the fiber class on each surface is the pullback of twice a line. For rational surfaces, the Picard group of the product splits as the sum of the two Picard groups. Thus $\mathcal O(X)$ is the pullback of $\mathcal O(2,2)$. Its defining section also descends, since $(\pi_i)_*\mathcal O_{S_i}=\mathcal O_{\mathbb P^2}$. We write the descended equation as $Q_a(z)=z^{\mathsf T}A(a)z$, with $a$ in the first plane and $z$ in the second.

After interchanging the factors, assume $S_2\ne\mathbb P^2$. Let $o$ be the first center in a sequence of blowups giving $\pi_2$. There is a point $a$ in the first plane satisfying
\[
 Q_a(o)=0,\qquad \det A(a)=0.
\]
Indeed, the two homogeneous forms have respective degrees two and six. If both are nonzero, their plane curves meet by B\'ezout's theorem. If either vanishes identically, the other still has a zero. Over any point of $\pi_1^{-1}(a)$ the fiber is defined on $S_2$ by the pullback of $Q_a$.

A zero quadratic form gives a surface fiber. Otherwise $Q_a$ defines a singular conic passing through $o$. For two distinct lines, the first blowup gives both strict transforms and the exceptional curve, even if only one line contains $o$. For a doubled line it gives the doubled strict transform and the doubled exceptional curve. Under subsequent blowups these components persist with their multiplicities. None of the resulting fibers is a plane conic. Thus at least one of the two general fiber classes has odd intersection with a divisor on the surface containing it.

Consequently, at least one projection $X\to S_i$ has a rational section. Its generic conic is a projective line over the function field of the base, so $X$ is birational to $S_i\times\mathbb P^1$ and is rational.
\end{proof}

In particular, $X$ is algebraically elliptic by Banecki's theorem \cite{Ban}.

\medskip

\par\noindent\textit{Remark.}
The embedding hypothesis in Theorem~\ref{thm:rational-surfaces} is essential for the rationality conclusion. Let
$q:X\to Y=(\mathbb P^1)^3$ be the double cover branched along a very general
smooth divisor of multidegree $(2,2,2)$. Then $X$ is smooth and
$-K_X=q^*\mathcal O_Y(1,1,1)$ is ample, so $X$ is rationally connected.
Composing $q$ with either projection to the factors $(1,2)$ or $(1,3)$ gives a
conic bundle over $\mathbb F_0$. Indeed, after fixing two coordinates,
the fiber is the plane conic $w^2=Q(u,v)$, where $Q$ is the binary quadratic
form obtained by restricting the branch equation to the remaining factor.

The product of these two projections factors through $q$ followed by the
embedding $(x,y,z)\mapsto((x,y),(x,z))$. Thus the map to the product of the
bases has degree two onto its image. A general conic of the first projection
maps with degree two onto a ruling line on the second base, so every divisor
pulled back from that base has even degree on the conic.
Hassett and Tschinkel prove that a very general such $X$ fails to be stably
rational \cite[\S7.4, following Lemma~8]{HTF}. For the algebraic ellipticity of these double covers,
see \cite{DKZ26b}.

\end{document}